\documentclass[11pt]{amsart}
\usepackage{geometry}                
\usepackage{amsfonts, amscd, amssymb, amsthm, amsmath}
\usepackage{pdfsync}
\usepackage{enumerate}
\usepackage{multirow}
\usepackage{float}

\theoremstyle{plain}
	\newtheorem{thm}{Theorem}[section]
	
	\newtheorem{prop}[thm]{Proposition}
	\newtheorem{cor}[thm]{Corollary}
\theoremstyle{definition}
	\newtheorem*{defn}{Definition}

\theoremstyle{example}

\usepackage{tikz}
\usepackage[pdftex,bookmarks]{hyperref}

\tikzstyle{V}=[draw, fill =black, circle, inner sep=0pt, minimum size=3pt]

\newcounter{r}
\newcommand\Part[1]{
        \setcounter{r}{1}
	 \foreach \x in {#1}{
 	{\ifnum\value{r}=1
		\draw (0,\value{r}-1)--(\x,\value{r}-1); 
		\fi}
	\draw (0,\value{r}) to (\x,\value{r});
   	\foreach \y in {0, ..., \x} {\draw (\y,\value{r})--(\y,\value{r}-1);}
	\addtocounter{r}{1}
 }}
 \def\PartUNIT{.175}
\newcommand{\PART}[1]{
\begin{matrix}
\begin{tikzpicture}[xscale=\PartUNIT, yscale=-\PartUNIT] 
	\Part{#1}
\end{tikzpicture}
\end{matrix}
}

\usepackage{color}

\def\cA{\mathcal{A}}\def\cC{\mathcal{C}}\def\cK{\mathcal{K}}\def\cU{\mathcal{U}}

  \def\CC{\mathbb{C}}   \def\FF{\mathbb{F}}        \def\NN{\mathbb{N}}   \def\QQ{\mathbb{Q}}         \def\ZZ{\mathbb{Z}}

\def\GL{\mathrm{GL}}

\def\<{\langle} \def\>{\rangle}

\def\End{\mathrm{End}}

\def\FI{\mathrm{FI}}

\title{A Survey of Representation Stability Theory}
\date{\today}                                           
\author{Anastasia Khomenko and Dhaniram Kesari}
\begin{document}
\maketitle
\setcounter{tocdepth}{1}
\tableofcontents

\begin{abstract}
In this survey article we summarize the current state of research in representation stability theory. We look at three different, yet related, approaches, using (1) the category of FI modules, (2) Schur-Weyl duality, and (3) finitely-generated modules over certain infinite dimensional vector spaces. The main example is the stability of representations of the symmetric group, though there have also been some notable generalizations of representation stability to other groups. This work summarizes the research that both authors engaged in over the course of the summer.
\end{abstract}

\section{Introduction}

Representation stability theory is the study of stable properties of representations of sequences of abstract algebraic structures, e.g. groups, rings, algebras, etc. Given the recent wave of new literature, we endeavor to describe the current state of research, collecting the main results, examples, and questions in the field.

Recall that a \textit{representation} of an algebra $A$ is a homomorphism from $A$ to the space of linear maps on a vector space $V$,
$$\rho: A \rightarrow \mathrm{End}(V).$$
In a sense, understanding the representation theory of an algebra is equivalent to understanding the algebra itself. But by focusing on the representations, often beautiful patterns are revealed in otherwise intractable algebraic structures. Representation stability theory focuses on patterns that stabilize for infinite families of algebras, such as the group algebras for the symmetric groups $S_n$, as $n$ varies. In \textit{combinatorial representation theory}, the representations are encoded using combinatorial tools like partitions and tableaux. Rules for combining representations are then given in terms of combinatorial data as well. Many of the favorite examples in representation stability are also favorite examples in combinatorial representation theory, like for the symmetric group, but also other Weyl groups, configuration spaces, and diagram algebras. For more details on combinatorial representation theory see (\cite{BR}).

One major surge in the field of representation stability theory came in the work of Church and Farb in \cite{CF}, where they defined representation stability in terms of consistent sequences to answer a number of questions about homological stability, cohomology of groups, Lie algebras and their homology, as well as combinatorics. They used their results to discuss not only representation theory, but a variety of problems including counting problems and finite group theory. Church, Ellenberg, and Farb then introduced the theory of $\FI$-modules in \cite{CEFa} to help explain representation stability of the symmetric group. In \cite{CEFb}, the authors transformed questions of stability into questions of finite generation. We outline the study of $\FI$-modules and their generalizations in Section \ref{FI}. A list of finitely generated $\FI$-modules can be found in Appendix \ref{app:A.1}.

Around the same time, Sam and Snowden began their own study of representation stability, which we discuss in Section \ref{TCA}. They took a similar approach to Church and Farb in \cite{SSa} and studied the algebraic structure of a category that they further developed in \cite{SSc}. They answered questions about how modules over twisted commutative algebras and their tensor products decompose. In \cite{SSa}, Sam and Snowden draw on Schur-Weyl duality in order to study the symmetric group via the downwards partition algebra. In \cite{SSc} the authors point out that the category they developed is equivalent to the category of $\FI$-modules in characteristic 0.
 
Bowman, De Visscher, and Orellana tackled the problem in \cite{BDO} by studying the stability of Kronecker coefficients. These are the decomposition numbers for the tensor product of symmetric group modules. As we explore in Section \ref{Kronecker}, Bowman, De Visscher, and Orellana utilize the Schur-Weyl duality between the symmetric group, $S_{n}$, and partition algebra, $P_{k}(n)$, transferring questions of stability for varying $S_{n}$'s to invariance of $P_{k}(n)$ with respect to its parameter $n$.

More classically, Murnaghan's Theorem describes certain stability of Kronecker coefficients explicitly. In \cite{Ma} and \cite{Mb}, Murnaghan studied the coefficients of the decomposition of the tensor product of two $S_n$-modules at length. Church, Farb, and Ellenberg re-proved Murnaghan's Theorem using their theory of $\FI$-modules in \cite{CEFb} and Wilson proved analogues of the theorem for the other Weyl groups of classical Lie type \cite{Wil}. Similarly, Bowman, DeVisscher, and Orellana were able to give yet another proof using the representation theory of the partition algebra. The Kronecker problem, that of giving positive combinatorial formulas for Kronecker coefficients, is one of the big open problems of the last century, so finding new approaches is of great interest to the community.

Several generalizations have also arisen in short succession. As mentioned above, Wilson generalized the theory of $\FI$-modules by replacing the symmetric group by other Weyl groups of classical Lie type in \cite{Wil} and \cite {W}. Putman and Sam studied categories in \cite{PS} that are analogous to $\FI$-modules by replacing the action of the symmetric group with other finite linear groups. Prior to the work of Church, Ellenberg, and Farb, Putman developed the notion of central stability in \cite{P}, which was shown to imply representation stability. Several researchers focused on exploring the cohomology of specific spaces using these new tools by considering actions of varying symmetric groups. In \cite{JR}, Rolland used the theory of $\FI$-modules to study the cohomology of pure mapping classes. In \cite{HR}, Hersh and Reiner explored the cohomology of configuration spaces. Nagpal discussed the cohomology of modular representations of symmetric groups in \cite{N} using the theory of $\FI$-modules.\\ 

\paragraph{\bf Acknowledgments.}

We would like to thank the Mathematics department at The City College of New York for their support through the Dr.\ Barnett and Jean Hollander Rich Mathematics Scholarships. We would also like to thank Chris Bowman, Jordan Ellenberg, and Steven Sam for helpful conversations and suggestions. We thank our research mentor Professor Zajj Daugherty for her guidance, unbridaled clarity, and endless patience.

\section{The $\FI$-module approach}\label{FI}

The development of $\FI$-modules was precipitated by the work of Church and Farb \cite{CF}, in which the authors introduced the notion of consistent sequences and a suitable definition of representation stability to answer a number of questions about homological stability, cohomology of groups, Lie algebras and their homology, as well as combinatorics. The following definitions are given in terms of a family of groups $\{G_n\}_{n \in \NN}$ with a fixed set of maps 
$$G_0 \hookrightarrow G_1 \hookrightarrow G_2 \hookrightarrow \cdots,$$
though they are largely motivated by the example where $G_n = S_n$ is the symmetric group, with the standard inclusion $S_n \hookrightarrow S_{n+1}$. 

\begin{defn} (\cite [Defn 1.1] {CF}) 
Let $V_n$ be a sequence of $G_n$-representations, equipped with linear maps $\phi_n: V_n \rightarrow V_{n+1}$, such that $\phi_n$ commutes with $g$, where $g \in G_n$ acts on $V_{n+1}$ by its image under the inclusion $G_n \hookrightarrow G_{n+1}$. Such a sequence is called \textit{consistent}. 

Now let $\left\{ V_n \right\}$ be a consistent sequence of $G_n$ representations. The sequence $\left \{ V_n \right\}$ is \textit{representation stable} if, for sufficiently large $n$, each of the following conditions hold.

\begin{enumerate}

\item{Injectivity:} The maps $\phi : V_n \rightarrow V_{n+1}$ are injective.
\item{Surjectivity:} The span of the $G_{n+1}$-orbit $of \phi _n(V_n)$ equals all of $V_{n+1}$.
\item{Multiplicities:} Decompose $V_n$ into irreducible $G_n$-representations as

$$V_n = \bigoplus_{\lambda} c_{\lambda ,n}V(\lambda)_n,$$

with multiplicities $0 \leq c_{\lambda ,n} < \infty .$

\end{enumerate}

\end{defn}

The theory of $\FI$-modules, defined as follows, was first presented by Church, Ellenberg, and Farb in \cite{CEFa} to encode sequences of $S_n$-representations with linear maps that preserve group structure. 

\begin{defn} 
An \textit{$\FI$-module} over a commutative ring $R$ is a functor,
$$V: \mathrm{FI} \rightarrow R\text{-Mod},$$ 

\noindent from the category whose objects are finite sets and whose morphisms are injective maps, to the category of modules over $R$.

\end{defn}

Note that the notion of $\FI$-modules is more restrictive than that of consistent sequences presented by Church and Farb in \cite{CF}. Specifically, if we have a consistent sequence of $S_n$ representation, we can transform it into an $\FI$-module if for all $m \leq n$: given two permutations $\sigma_1, \sigma_2 \in S_n$ and an element $v \in V_n$, we have
$$ \sigma_1|_{\left\{1,2,...,m \right\}} = \sigma_2|_{\left\{1,2,...,m \right\}} \Rightarrow \sigma_1(v) = \sigma_2(v).$$

\subsection{Finite generation of $\FI$-modules}

\begin{defn}
An $\FI$-module is \textit{finitely generated} if there is a finite set $S$ of elements in $\coprod_i V_i$ so that no proper sub-$\FI$-module of $V$ contains $S$.

\end{defn}

Several examples appear in \ref{app:A.1}. \\

As we will see in Theorem 2.2, representation stability properties can be transferred to finite generation properties of $\FI$-modules, so it was important to establish the relation between those two concepts. In order to do so, Church, Ellenberg, and Farb first proved the Noetherian property of certain $\FI$-modules and then used it to link finite generation of those $\FI$-modules to representation stability of the corresponding $S_n$ representations.

\begin{thm} (\cite [Thm 2.60] {CEFa}) \textbf{($\FI$-modules are Noetherian)} 
Let $R$ be a Noetherian ring containing $\mathbb{Q}$. The category of $\FI$-modules over $R$ is Noetherian, that is, any sub-$\FI$-module of a finitely generated $\FI$-module is finitely generated.

\end{thm}

The relation between finite generation of $\FI$-modules and representation stability is given by the following theorem.

\begin{thm} (\cite [Thm 1.14] {CEFa}) \textbf{(Finite generation vs Representation stability)} An $\FI$-module $V$ over a field of characteristic 0 is finitely generated if and only if the sequence $\{ V_n \}$ of $S_n$ representations is uniformly representation stable (see \cite[Defn 2.3]{CF}), and each $V_n$ is finite-dimensional. In particular, for any finitely-generated $\FI$-module $V$, we have for sufficiently large $n$ a decomposition
$$ V_n \simeq \bigoplus c_{\lambda}V(\lambda)_n.$$

\end{thm}

\subsection{Generalizations: $\FI$-modules of other Lie types.}

Later, the theory of $\FI$-modules was extended by Wilson in \cite{Wil} and \cite{W} to include the representation stability of other classical Weyl groups. Such $\FI$-modules are denoted as $FI_W$-modules in general (where $W$ is a sequence of either $S_n$ or $BC_n$ or $D_n$ groups) and $FI_A$, $FI_{BC}$ and $FI_D$ in particular. For background on Weyl groups of classical type see \cite[Section 2]{Wil} 

Just as in the symmetric group (type A) case, Wilson showed that finite generation of $FI_{BC}$- and $FI_D$-modules is equivalent to representation stability (as defined in \cite{CF}) of $BC_n$ and $D_n$ (\cite{Wil}). 
Wilson also proved that the restriction of a finitely generated $FI_W$-module preserves finite generation.

\begin{defn}(\cite [Defn 1.3] {Wil}
An $FI_W$-module is finitely \textit{generated in degree $\leq d$}, if it has a finite generating set $\left\{ v_i \right\}$ with $v_i \in V_{m_i}$, for some finite sequence of integers  $\{ m_i \}$, where $m_i \leq d$ for each $i$.
\end{defn}

Let $W$, $\overline{W}$ be families of Weyl groups such that $W_n \subseteq  \overline{W}_n$.

\begin{defn}(\cite [Defn 3.21] {Wil}
Given a family of inclusions $W_n \hookrightarrow \overline{W}_n$, any $\FI_{\overline{W}}$-module $V$ inherits the structure of an $FI_W$-module by restricting the functor $V$ to the subcategory $\FI_{W}$ in $\FI_{\overline{W}}$. We call this restriction $Res^{\overline{W}}_{W}V$ the restriction of V to $FI_W$.

\end{defn}

\begin{prop}(\cite[Proposition 3.22] {Wil}).\textbf{ (Restriction preserves finite generation)}
For each family of Weyl groups $W \subset \overline{W}$, the restriction $Res^{\overline{W}}_{W}V$ of a finitely generated $FI_{\overline{W}}$-module V is finitely
generated as an $FI_{W}$-module. Specifically, we have the following.\\
1. Given an $FI_{\emph{BC}}$-module $V$ finitely generated in degree $\leq m$, $Res^{BC}_{A}V$ is finitely generated as an $FI_{\emph{A}}$-module in degree $\leq m$.\\
2. Given an $FI_{\emph{BC}}$-module $V$ finitely generated in degree $\leq m$, $Res^{BC}_{D}V$ is finitely generated as an $FI_{\emph{D}}$-module in degree $\leq m$.\\
3. Given an $FI_{\emph{D}}$-module $V$ finitely generated in degree $\leq m$, $Res^{D}_{A}V$ is finitely generated as an $FI_{\emph{A}}$-module in degree $\leq m+1$.
\end{prop}

 She also showed that induction respects generation and relation degree. The natural embeddings $S_n \hookrightarrow D_n \hookrightarrow BC_n$ provide us with the inclusions of categories,

$$FI_A \hookrightarrow FI_D \hookrightarrow FI_{BC}.$$

\begin{defn}

An FI module $V$ has $\textit{relation degree} \leq r$ if there exists a surjection
$$\bigoplus_{i=0}^d M(W_i) \twoheadrightarrow V,$$

\noindent whose kernel $K$ is generated in degree $\leq r$.

\end{defn}

\begin{thm} (\cite [Cor 3.28] {Wil}) \textbf{(Induction preserves finite generation).}
Suppose $V$ is a finitely generated $FI_W$-module with degree of generation $\leq g$ and relation degree $\leq r$. Then $Ind^{\overline{W}}_WV$ is also finitely generated, with generation degree $\leq g$ and the relation degree $\leq r$.

\end{thm}

The above result provided the sufficient groundwork to show that the dimensions of finitely generated $\FI_{BC}$- and $\FI_D$-modules over arbitrary fields are eventually polynomial, and to prove the Noetherian property for $\FI_{BC}$- and $\FI_D$-modules. The latter extends the theorem about finite generation and representation stability proved by Church, Ellenberg, and Farb in \cite{CEFa} to include sequences of $BC_n$ and $D_n$ groups. This result, as well as specific stable ranges, is presented in the two following theorems.

\begin{defn}  
Let $V$ be an $\FI$-module over a field of characteristic 0. We say that $V$ has \textit{weight} $\leq d$ if for every $n \geq 0$ and every irreducible constituent $V(\lambda)_n$ of $V_n$, we have $|\lambda| \leq d$.
\end{defn}

\begin{thm} 

(\cite [Thm 4.27] {Wil}) \textbf{(Finitely generated $\FI_W$-modules are uniformly representation stable).} Suppose that $\FF$ is a field of characteristic 0, and $W_n$ is $S_n$, $D_n$, or $BC_n$. Let $V$ be a finitely generated $\FI$-module. Take $d$ to be an upper bound on the weight of $V$, $g$ an upper bound on its degree of generation, and $r$ an upper bound on its relation degree. Then ${V_n}$ is uniformly representation stable with respect to the maps induced by the natural inclusions
$$I_n: n \hookrightarrow (n+1),$$

\noindent stabilizing once $n \geq \text{max}(g,r)+d$; when $W_n$ is $D_n$ and $d=0$, we need the additional condition that $n \geq g+1$.

\end{thm}

\begin{thm}
(\cite [Thm 4.28] {Wil}) \textbf{(Uniformly representation stable $\FI_W$-modules are finitely generated).} Suppose conversely that $V$ is an $\FI_W$-module, and that $\{V_n, (I_n)_*\}$ is uniformly representation stable for $n \geq N$, where $(I_n)_*$ is the following inclusion,
$(I_n)_*: V_n \hookrightarrow V_{n+1}.$ Then V is finitely generated in degree $\leq N.$

\end{thm}

After obtaining this result, in her next paper \cite {W}, Wilson provided a specific stable range of sequences of $W_n$-representations, such as cohomology groups of pure string motion groups $P \Sigma_n$ and complexified complements $M_W(n)$ of the reflecting hyperplanes of $W_n$, that are associated with $\FI_W$-modules. She also proved that the characters of these representations are generated by polynomials. 

\begin{thm}(\cite [Thm 1.1] {W})\textbf{(Polynomiality of characters).}
 Let $W_n$ denote $S_n, B_n$, or $D_n$. Let $X_n$ denote either the sequence of pure string motion groups $P \Sigma_n$ or the complexified complements $M_W(n)$ of the reflecting hyperplanes of $W_n$. Then, in each degree $i$ the sequence of cohomology groups of such sequences $\{ H^i(X_n; \mathbb{Q}) \}$ is an $\FI_W$-module finitely generated in degree $\leq 2i$. Then,

\begin{enumerate}
\item the sequence $ \{ H^i(X_n; \mathbb{Q}) \}$ is uniformly representation stable, stabilizing for $n \geq 4i$; and\\
\item for all values of $n$, the characters of $ \{ H^i(X_n; \mathbb{Q}) \}$ are given by a unique character polynomial of degree at most $2i$. 
\end{enumerate}

\end{thm}

\noindent Additionally, Wilson proved that the characters $BC_n$-representations in the cohomology rings of pure string motion groups $H^\ast(P \Sigma_\bullet; R)$ are given by a polynomial whose degree is bounded from above.

\begin{cor} (\cite [Cor 5.6] {W})
Fix an integer $i \geq 0$. The character of the sequence of $B_n$-representations $\{ H^i(P \Sigma_n; \mathbb{Q}) \} _n$ are given by a unique character polynomial of degree $\leq 2i$ for all values of $n$. 

\end{cor}

\subsection{Other techniques resulting from the theory of $\FI$-modules}

The development of $\FI$-modules gave rise to other techniques used to prove representation stability of the symmetric group $S_n$. One such method was presented in Hersh and Reiner's paper, in which they study the stability of the $S_n$-representations given by cohomology of configuration spaces in $\mathbb{R}^d$. 

The first theorem of the paper extended the work of Church (\cite[Thm 1] {C}) by providing the specific $n$ for which the $S_n$-representations in $H^i(\mathrm{Conf}(n, \mathbb{R}^d))$ stabilize sharply. 

\begin{thm} (\cite[Thm 1.1] {HR})
Fix integers $d\geq 2$ and $i \geq 1$. Then $H^i(\mathrm{Conf}(n, \mathbb{R}^{d}))$ vanishes unless $d-1$ divides $i$, in which case it stabilizes sharply at 
$$ n = \begin{cases}

3 \frac{i}{d-1} \text{ for } d\geq 3 \text{ odd,}\\
3 \frac{i}{d-1} +1 \text{ for } d\geq 2 \text{ even. }

\end{cases}$$

\end{thm}

This result was achieved by recasting the cohomology groups $H^i(\mathrm{Conf}(n, \mathbb{R}^d))$ in terms of higher Lie characters $\mathrm{Lie}_{\lambda}$ when $d$ is odd, and the Whitney homology of the lattice $\Pi_n$ of set partitions of $\left\{ 1,2,...,n \right\}$, when $d$ is even. This was done the following way
$$H^i(\mathrm{Conf}(n, \mathbb{R^d})) \cong \begin{cases}
\bigoplus_{m=i+1}^{2i}M_n(\widehat{Lie}^i_m)   \text{ for } d \text{ odd,}\\
\bigoplus_{m=i+1}^{2i}M_n(\widehat{W}^i_m) \text{ for } d \text{ even,}
\end{cases}$$

\noindent where $\widehat{Lie}^i_m$ and $\widehat{W}^i_m$ are certain subrepresentations of higher Lie characters (i.e. characters that result from the Poincare-Birkhoff-Witt basis of the free associative algebra and which sum up to regular Lie characters of $S_n$ representations) (\cite[section 2.6] {HR}) and Whitney homology of the lattice $\Pi_n$ of set partitions of   $\left\{ 1,2,...,n \right\}$ (\cite[section 2.4] {HR}). \bigskip

Although the irreducible decompositions of $\widehat{Lie}^i_m$ and $\widehat{W}^i_m$ are still not known in general, the authors provided the irreducible decompositions of sums of $\widehat{Lie}^i_m$ and $\widehat{W}^i_m$. 

The sums of irreducible representations $\widehat{Lie}_n^i$ and $\widehat{W}^i_n$ are respectively denoted as $\widehat{Lie}_n$, and $\widehat{W}_n$. 

\begin{thm} (\cite[Thm 1.3] {HR})
These sums have the irreducible decompositions,
$$\widehat{Lie}_n = \Sigma_Q \chi^{shape(Q)} \quad \mathrm{and} \qquad \widehat{W}_n = \Sigma_Q \chi^{shape(Q)}, $$

\noindent in which the sums range over the set of desarrangement tableaux (standard tableaux $Q$ with even first ascent), and Whitney-generating tableaux $Q$ of size $n$ (\cite[Defn 6.3] {HR}).\end{thm}

The latest development of the study of $\FI$-modules brings us the result obtained by Li in (\cite[Thm 1.17] {LL}) which provides us with information about an upper bound of homological degree of finitely generated $\FI$-modules,

\begin{defn}

The $i^{th}$ homological degree of V and is defined as follows:
$$hd_i(V) = \mathrm{sup}\{ m \in \ZZ_{+} ~|~ \mathrm{the~value~of~}H_i(V) \mathrm{~on }~ m \mathrm{~is~not~} 0 \}$$

\end{defn}

\begin{thm} (\cite[Thm 1.17] {LL})

Let $\FF$ be a field of characteristic 0, and let V be a finitely generated $\FI$-module. Then for $i \geq 1$, 

$$\mathrm{hd}_i(V) \leq \mathrm{max} \{2\mathrm{gd}(V) - 1, \mathrm{td}(V)\} + \mathrm{i},$$

where $\mathrm{gd}(V)$ is the homolological degree $hd_0(V)$, and $\mathrm{td}(V)$ is the torsion degree of $V$. 
\end{thm}

In \cite{R}, Ramos extended the results in \cite[Thm B] {CEFN} to $FI_G$-modules and  provided specific bounds for the stable range of these modules.

\begin{defn}
A $G$-map between sets $R$ and $S$ is a pair ($f$, $\rho$), where $f: R \rightarrow S$ and 
$\rho: R \rightarrow G$.

\end{defn}

\begin{defn}
Let $\mathrm{FA}_G$ be the category whose objects are finite sets and whose morphisms are $G$-maps. 
Then define $\FI_G$ to be the subcategory where the function $f$ is injective. 

\end{defn}

\begin{thm}(\cite[Thm D] {R}) \textbf{(Stable range of $\FI_G$-modules)}
Let $G$ be a finite group, and let $V$ be a finitely generated $\FI_G$-module over a field $\FF$. Then the stable range of $V$ $\geq r + min\{ r, d\}$, where $r$ is the relation degree of $V$ and $d$ is the degree of the $0^{th}$ homology $H_0(V)$.

\end{thm} 

\subsection{Extending the theory of $\FI$-modules} 

In \cite{PS}, Putman and Sam constructed analogues of $\FI$-modules, where the role of the symmetric groups is played by the general linear groups and the symplectic groups over finite rings. The Asymptotic Structure Theorem established an analogue of the Noetherian property of $\FI$-modules. The results of \cite{PS} demonstrated the importance of Noetherianity over finitely generated rings.

\begin{defn}{\cite[Central stability.]{P}} For each $n$, let $V_{n}$ be a representation of $S_{n}$ over a field $\FF$, and let $\phi_{n}: V_{n} \rightarrow V_{n+1}$ be a linear map which is $S_{n}$-equivariant. We will call the sequence
$$ V_{1} \xrightarrow{\phi_{1}} V_{2} \xrightarrow{\phi_{2}} V_{3} \xrightarrow{\phi_{3}} V_{4}\xrightarrow{\phi_{4}} \hdots
$$
a \emph{coherent sequence} of representations of the symmetric group if it forms a sequence of equivariant linear maps. The central stabilization of $\phi_{n-1}$, denoted $\cC(V_{n-1} \xrightarrow{\phi_{n-1}} V_{n})$, is the $S_{n+1}$-representation which is the largest quotient of $Ind^{Sn+1}_{Sn}V_{n}$ such that $(n,n+1)$ acts trivially on the image of $V_{n-1}$.We will say that our coherent sequence is \emph{centrally stable} starting at $N \geq 2$ if for all $n \geq N$, we have $V_{n+1} = \cC(V_{n-1} \xrightarrow{\phi_{n-1}} V_{n})$ and $\phi_{n}$ is the natural map $V_{n} \rightarrow \cC(V_{n-1} \xrightarrow{\phi_{n-1}} V_{n})$.
\end{defn}

\begin{defn}{\cite{PS}} Let $R$ be a ring. A linear map between free $R$-modules is \emph{splittable} if its cokernel is free. Let $VI(R)$ be the category whose objects are finite-rank free $R$-modules and whose morphisms are splittable injections. For $V \in VI(R)$, the monoid of $VI(R)$-endomorphisms of $V$ is $GL(V)$. Thus if $M$ is a $VI(R)$-module, then $M_{V}$ is a representation of $GL(V)$ for all $V \in VI(R)$.
\end{defn}

\begin{thm}{\cite[Theorem A]{PS}} Let $R$ be a finite ring. The category of $VI(R)$-modules is Noetherian.
\end{thm}

Putman and Sam defined two other categories, namely the categories $VIC$ and $SI$. Interested readers should visit {\cite[Section 1]{PS}} for the definitions of these categories. They also proved Noetherianity for the category $OVIC(R)$ in {\cite[Section 2.3]{PS}}. Ultimately, the authors showed that both of these are Noetherian. The authors generalized their results by proving the Asymptotic Structure Theorem. It should be noted that Church, Farb, Ellenberg, and Nagpal studied the asymptotic structure of $\FI$-modules in \cite{CEFN}. The Asymptotic Structure Theorem Putman and Sam proved is able to provide a broadened approach to many of the ideas explored in the study of $FI$-modules.

\begin{defn}{\cite{PS}} A \emph{weak complemented category} is a monoidal category
$(A,\circledast)$ satisfying the following.
\begin{enumerate}[(1)]
\item Every morphism in $A$ is a monomorphism. Thus for all morphisms $f: V \rightarrow V'$, it makes
sense to talk about the subobject $f(V)$ of $V'$.
\item The identity object 1 of $(A,\circledast)$ is initial, i.e.\ for $V$, $V' \in A$ there exist natural morphisms
$V \rightarrow V \circledast V'$ and $V' \rightarrow V \circledast V'$, and thus the compositions
\end{enumerate}
$$
V \xrightarrow{\cong} V \circledast 1 \rightarrow V \circledast V'
\text{ $\qquad$ and $\qquad$} 
V' \xrightarrow{\cong} 1 \circledast V' \rightarrow V \circledast V'.
$$
We will call these the \emph{canonical morphisms}. Note the following.
\begin{enumerate}[(1)]
\item For $V$, $V'$, $W \in A$, the map 
$$Hom_{A}(V \circledast V', W) \rightarrow Hom_{A}(V, W) \times Hom_{A}(V',W),$$
obtained by composing morphisms with the canonical morphisms, is an injection.\\
\item Every subobject $C$ of an object $V$ has a unique complement. So, a subobject $D$ of $V$ such that there is an isomorphism $C \circledast D \xrightarrow{\cong} V$, where the compositions $C \rightarrow C \circledast D \xrightarrow{\cong} V$ and $D \rightarrow C \circledast D \xrightarrow{\cong} V$ of the isomorphism with the canonical morphisms are the inclusion morphisms.
\end{enumerate}
\end{defn}

\begin{defn}{\cite{PS}} A \emph{complemented category} is a weak complemented category $(A,\circledast)$ whose monoidal structure
is equipped with a symmetry.
\end{defn}

\begin{defn}{\cite{PS}} If $(A,\circledast)$ is a monoidal category, then a generator for $A$ is an object $X$ of $A$ such that all objects $V$ of $A$ are isomorphic to $X_{i}$ for some unique $i \geq 0$. We will call $i$ the $X$-rank of $V$.
\end{defn}
The following theorem is a major result of \cite{PS} since it generalizes the criteria for Noetherianity.

\begin{thm}{\cite[Theorem E]{PS}} Let $(A,\circledast)$ be a complemented category with generator $X$. Assume that the category of $A$-modules is Noetherian, and let $M$ be a finitely generated $A$-module. Then the following hold. For $N \geq 0$, let $A^{N}$ denote the full subcategory of $A$ spanned by elements whose $X$-ranks are at most N.\\
$(1)$ (Injective representation stability) If $f : V \rightarrow W$ is an $A$-morphism, then the homomorphism
$M_{f} : M_{V} \rightarrow M_{W}$ is injective when the $X$-rank is sufficiently large. \\
$(2)$ (Surjective representation stability) If $f : V \rightarrow W$ is an $A$-morphism, then the orbit under $Aut_{A}(W)$ of the image of $M_{f} : M_{V} \rightarrow M_{W}$ spans $M_{W}$ when the $X$-rank is sufficiently large. \\
$(3)$ (Central stability) For $N \gg 0$, the functor $M$ is the left Kan extension to FI of the restriction of $M$ to $A^{N}$.
\end{thm}

They are able to prove a version of this theorem that is specific to $\FI$-modules based on the work in \cite{CEFN} where $A$ is FI and the generator $X$ is $\{1\}$ of the prior theorem.

As the authors noted, their approach to proving the Noetherian property works specifically for categories over finite rings.

\begin{thm}{\cite[Theorem M]{PS}}
For $C \in \{VI(\ZZ), VIC(\ZZ), SI(\ZZ)\}$, the category of $C$-modules is not Noetherian.
\end{thm}

\section{The centralizer algebra approach}\label{Kronecker}

One particularly elusive question in representation theory addresses the Kronecker coefficients, the decomposition numbers for the tensor product of two symmetric group modules. Namely, if we write $S^\lambda$ for the irreducible $S_n$ representation indexed by the partition $\lambda$, then in
$$ S^{\lambda} \otimes S^{\mu} = \bigoplus_{\nu \vdash n} c_{\lambda, \mu}^{\nu} S^{\nu} \text{, } \qquad c_{\lambda, \mu}^{\nu} \in \ZZ_{\geq 0},$$ 
the coefficients $c_{\lambda, \mu}^{\nu}$ are \emph{Kronecker coefficients}. One major contribution to this open question is Murnaghan's Theorem, describing their behavior for varying symmetric groups. In Bowman, De Visscher, and Orellana's study of representation stability theory in \cite{BDO}, the authors focused on these coefficients. In \cite{SSf} Sam and Snowden were able to prove Stembridge's Conjecture, a conjecture that generalizes Murnaghan's Theorem. The Kronecker problem, finding a combinatorially positive formula for the Kronecker coefficients remains frustratingly open. In \cite{BDO}, key insights are provided.

\subsection{Schur-Weyl Duality} 

Schur-Weyl duality offers a different lens through which we can study the symmetric group. The original result establishes the relationship between irreducible finite-dimensional representations of the general linear and symmetric groups. Schur found that the general linear group $\GL_n(\CC) = \GL_n$ and the symmetric group $S_k$ have commuting actions on the $k$-fold tensor product $(\CC^n)^{\otimes k} = \CC^n \otimes \hdots \otimes \CC^n$, where the general linear group acts diagonally and the symmetric group acts by permuting factors. Extending these actions linearly to the respective group algebras, these actions fully centralize each other in $\End(V^{ \otimes k})$, which yields the multiplicity-free decomposition
$${(\CC^n)}^{\otimes k} = \bigoplus_{\lambda \in \Lambda} G^\lambda \otimes S^\lambda
\qquad \text{as a $(\CC \GL_n, \CC S_k)$-bimodule,}$$
where $\Lambda$ is the set of integer partitions of $k$ with no more than $n$ parts, $S^{\lambda}$ is an irreducible $S_n$ module, and $G^{\lambda}$ is an irreducible 
$\GL_n$ module. Thus we have a pairing general linear and symmetric group modules; namely, 
for each $k \in  \ZZ_{\geq 0}$ and partition $\lambda$ of $k$ of length less than $n$, there is a distinct finite dimensional irreducible representation of $\GL_n$.

\subsection{The Partition Algebra} 

Like $\CC S_k$ is the centralizer of the diagonal action of $\CC\GL_{n}$ on $(\CC^n)^{\otimes k}$, the partition algebra $P_{k}(n)$ is the centralizer algebra of the diagonal action of $\CC S_{n}$ on  $(\CC^n)^{\otimes k}$, where $\CC^n$ is now the $n$-dimensional permutation representation of $S_n$.  The corresponding $(P_{k}(n), \CC S_{n})$-bimodule decomposes as 
$$ (\CC^n)^{\otimes k} = \bigoplus_{\mu \in \Lambda} P^{\mu} \otimes S^{\overline{\mu}},$$
where $\Lambda$ is the set of integer partitions of $\{1, \hdots, k\}$ with first part no longer than $n/2$, $P^\mu$ is an irreducible $P_k(n)$-module, and for $\mu\in \Lambda$,  $\overline{\mu}$ is the partition $(n - |\mu|, \mu_1, \mu_2, \dots)$. For example, if 
$$\mu = \PART{4,2,1} \qquad \text{ and } \qquad n=15,
 \qquad \text{then} \qquad \bar{\mu} = \PART{8,4,2,1}.$$ 
 Bowman, DeVisscher, and Orellana rigorously formulate the representation theory of this bimodule in Section 2 of \cite{BDO}.

The partition algebra $P_k(n)$ also has a remarkable description as a diagram algebra. Namely, it is the algebra with basis indexed by \emph{diagrams} on $2k$ vertices: plane graphs with $k$ vertices on top, $k$ vertices on bottom, and edges connecting vertices to represent connected components corresponding to a set partition of the $2k$ vertices. For example, 
$$
d_1 = 
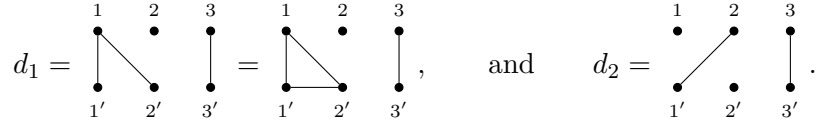
\begin{matrix}
	\begin{tikzpicture}[scale=.75]
	\foreach \x in {1,2,3}{
		\node[V, label=above:{\tiny$\x$}] (t\x) at (\x,1){};
		\node[V, label=below:{\tiny$\x'$}] (b\x) at (\x,0){};}
	\draw (t1) to (b1)  (b2) to (t1) (t3) to (b3);
	\end{tikzpicture}
\end{matrix}= \begin{matrix}
	\begin{tikzpicture}[scale=.75]
	\foreach \x in {1,2,3}{
		\node[V, label=above:{\tiny$\x$}] (t\x) at (\x,1){};
		\node[V, label=below:{\tiny$\x'$}] (b\x) at (\x,0){};}
	\draw (t1) to (b1) to (b2) to (t1) (t3) to (b3);
	\end{tikzpicture}
\end{matrix}, \qquad 
\text{and} \qquad
d_2 = \begin{matrix}
	\begin{tikzpicture}[scale=.75]
	\foreach \x in {1,2,3}{
		\node[V, label=above:{\tiny$\x$}] (t\x) at (\x,1){};
		\node[V, label=below:{\tiny$\x'$}] (b\x) at (\x,0){};}
	\draw (t2) to (b1) (b3) to (t3);
	\end{tikzpicture}
\end{matrix}.
$$
Composition of diagrams is performed as follows: place $d_1$ above $d_2$, connect the bottom vertices of $d_1$ with the top vertices of $d_2$ and remove all the connected components that are completely situated in the middle of the two diagrams. For example, if $d_1$ and $d_2$ are as above, then
$$
d_1 \circ d_2 = 
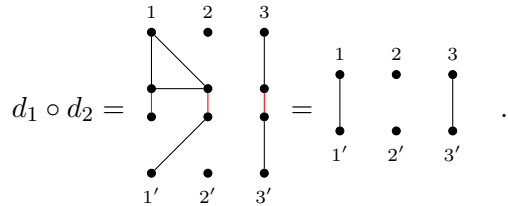
\begin{matrix}
	\begin{tikzpicture}[scale=.75]
	\foreach \x in {1,2,3}{
		\node[V, label=above:{\tiny$\x$}] (tt\x) at (\x,2.5){};
		\node[V] (tb\x) at (\x,1.5){};
		\node[V] (t\x) at (\x,1){};
		\node[V, label=below:{\tiny$\x'$}] (b\x) at (\x,0){};
		\draw[red] (t\x) to (tb\x);}
	\draw (tt1) to (tb1) to (tb2) to (tt1) (tt3) to (tb3);
	\draw (b1) to (t2) (b3) to (t3);
	\end{tikzpicture}
\end{matrix}
= \begin{matrix}
	\begin{tikzpicture}[scale=.75]
	\foreach \x in {1,2,3}{
		\node[V, label=above:{\tiny$\x$}] (t\x) at (\x,1){};
		\node[V, label=below:{\tiny$\x'$}] (b\x) at (\x,0){};}
	\draw (t1) to (b1) (t3) to (b3);
	\end{tikzpicture}
\end{matrix}\quad .
 $$
Then multiplication is given by $d_1d_2 = n^l(d_1 \circ d_2)$, where $l$ is the number of removed components. See, for example, \cite [Section 2] {HaRa}.

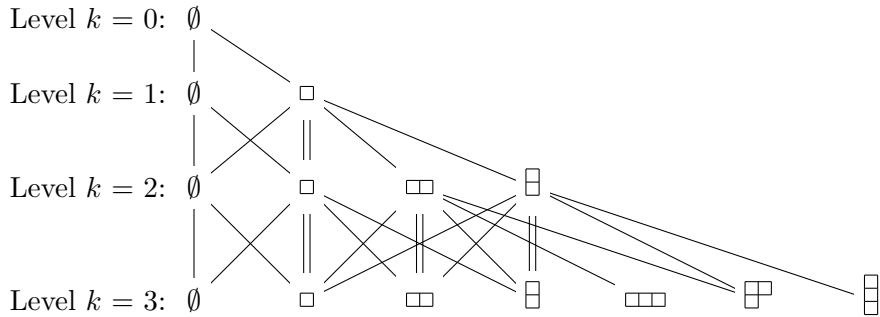
\begin{figure}[H]
\begin{tikzpicture}[yscale=-1, xscale=1.5]
	\coordinate (00) at (0,0);
	\coordinate (10) at (0,1);
	\coordinate (11) at (1,1);
	\foreach \x in {0, ..., 3}{\coordinate (2\x) at (\x,2.25);}
	\foreach \x in {0, ..., 6}{\coordinate (3\x) at (\x,3.75);}
\foreach \y in {0, ..., 3} { \node[left] at (\y0) {Level $k$ = \y:~~~.};}
\draw (00)--(10) (00)--(11);
\foreach \x in {0,1}{ \draw (10)--(2\x);}
	\foreach \x in {0, 2, 3}{ \draw (11)--(2\x);}
	\draw[double distance = 2pt] (11)--(21);
\foreach \x in {0,1}{ \draw (20)--(3\x);}
	\foreach \x in {1,2,3}{\draw[double distance = 2pt] (2\x)--(3\x);}
	\foreach \x in {0, 2, 3}{ \draw (21)--(3\x);}
	\foreach \x in {1, 3, 4, 5}{ \draw (22)--(3\x);}
	\foreach \x in {1, 2, 5, 6}{ \draw (23)--(3\x);}
\begin{scope}[every node/.style={fill=white}]
	\node at (00) {$\emptyset$};
	\node at (10) {$\emptyset$};
	\node at (11) {$\PART{1}$};
	\node at (20) {$\emptyset$};
	\node at (21) {$\PART{1}$};
	\node at (22) {$\PART{2}$};
	\node at (23) {$\PART{1,1}$};
	\node at (30) {$\emptyset$};
	\node at (31) {$\PART{1}$};
	\node at (32) {$\PART{2}$};
	\node at (33) {$\PART{1,1}$};
	\node at (34) {$\PART{3}$};
	\node at (35) {$\PART{2,1}$};
	\node at (36) {$\PART{1,1,1}$};
	\end{scope}
\end{tikzpicture}
\caption{\bf Bratteli diagram for $P_0(n) \subseteq P_1(n) \subseteq P_{2}(n) \subseteq \cdots$, levels 0--3.}
\label{fig:PkBrat}
\end{figure}

A useful tool for encoding the duality between the symmetric group and partition algebra is the Bratteli diagram. The first few levels of the Bratteli diagram for generic $n$ are presented in Figure~\ref{fig:PkBrat}. The paths from $\emptyset$ at the $0^{\text{th}}$ level to a point $\lambda$ on the $k^{\text{th}}$ level provide a basis for the simple $P_k(n)$-module, $P^\lambda_k(n)$, for $n\geq 2k$. Therefore, they record the multiplicity of $S^\lambda$ in a decomposition of $(\CC^n)^{\otimes k}$ as a $S_n$-module for $n\geq 2k$. 
In particular, for $n \geq 2k$, the representations of $P_k(n)$ are independent of $n$. 

%
%

The goal of Bowman, De Visscher, and Orellana in \cite{BDO} is to pass questions about the stability of the Kronecker coefficients to the partition algebra. They are thus able to offer an alternative approach to Murnaghan's Theorem stemming from stability of $P_{k}(n)$ with respect to varying the parameter $n$, $i.e.$ the representation theory of $P_{k}(n)$ is largely independent of $n$. The authors offered a procedure for how to calculate the Kronecker coefficients through the partition algebra, providing specific closed formulas for hook and two-part partitions.

\begin{thm}{\cite[Corollary 5.2]{BDO}}
Let $\lambda_{[n]}$, $\mu_{[n]}$, $\nu_{[n]}$ be partitions of $n$ with $| \lambda |$ = $r$, $|\mu|$ = $s$ and
$|\nu| = r + s - l$.\\
\indent $(i)$ Suppose $\nu_{[n]} = (n - k, k)$ is a two-part partition. Then we have
$$
g^{(n-k,k)}_{\lambda_{[n]}, \mu_{[n]}} =\sum_{\substack{l_{1}, l_{2} \\ l = l_{1} + l_{2}}} \sum_{\substack{\sigma \vdash l_{1} \\ \gamma \vdash l_{2}}} c^{\lambda}_{(r-l_{1}-l_{2}), \sigma , \gamma} c^{\mu}_{ \gamma , \sigma , (s-l_{1}-l_{2})}
$$
for all $n \geq$ min$\{ |\lambda|$ + $\mu_{1}$ + $k, |\mu| + \lambda_{1} + k \}$.\\
\indent $(ii)$ Suppose $\nu_{[n]} = (n - k, 1^{k})$ is a hook partition. Then we have
$$
g^{(n-k, 1^{k})}_{\lambda_{[n]}, \mu_{[n]}} =\sum_{\substack{l_{1}, l_{2} \\ l = l_{1} + 2l_{2}}} \sum_{\substack{\sigma \vdash l_{1} \\ \gamma \vdash l_{2}}} c^{\lambda}_{( 1^{r} - l_{1} - l_{2}), \sigma , \gamma} c^{\mu}_{ \gamma , \sigma ' , (1^{s}-l_{1}-l_{2})}
$$
for all $n \geq$ min$\{|\lambda| + |\mu| + 1, |\mu| + \lambda_{1} + k, |\lambda| + \mu_{1} + k \}$ and where $\sigma '$ denotes the transpose of $\sigma$.
\end{thm}

Just as Wilson proved analogues of Murnaghan's Theorem for the other Weyl groups, a natural question is to consider how this approach can be adapted to the other Weyl groups. It may also be interesting to use the Schur-Weyl duality between the symmetric group and partition algebra to study the representation theory of FI-modules through the lens of the partition algebra, since the partition algebra is largely $n$ invariant. In the following section, we will see Sam and Snowden make use of the $P_{k}(n), \CC S_{n}$-bimodule as well.

\section{The twisted commutative algebra approach}\label{TCA}

Sam and Snowden approached the problem of describing representation stability \cite{SSa} by analyzing modules over twisted commutative algebras (TCAs). They summarize their investigation of TCAs in \cite{SSb} and offer several equivalent definitions for TCAs and the modules over them. 
\begin{defn}{\cite[Definition 1]{SSb}} A \emph{twisted commutative algebra} is a rule which associates to each vector space $V$ a commutative ring $A(V)$, and to each linear map of vector spaces $V \rightarrow V ′$ a ring homomorphism $A(V) \rightarrow A(V′)$.
\end{defn}
Many TCAs satisfy the Noetherianity property making their representation theory tractable. In fact, Nagpal, Sam and Snowden have recently explored the Noetherian property over twisted commutative algebras in greater detail in \cite{SSe}, specifically the Noetherianity of  $\text{Sym(Sym}^{2}({\CC}^{\infty})).$ 

A particularly interesting TCA is $A=\text{Sym}(\CC\<1\>)$, the symmetric algebra over the vector space $\CC^{\infty}$. In \cite{SSa}, the authors note that the TCA Sym($\CC\<1\>$) is studied by Church, Farb, and Ellenberg as FI-modules. In \cite{CEFa}, the authors state that "FI-modules can be viewed as modules for the `exponential' twisted commutative algebra". This TCA is studied at length in \cite{SSa}. Some important structural results include the following.

\begin{thm}{\cite[Cor 4.2.6]{SSa}}
Projective $A$-modules are also injective.
\end{thm}
 
\begin{thm}{\cite[Thm 4.3.1]{SSa}}
Every object of $Mod_A$ has finite injective dimension.
\end{thm}

In \cite{SSc}, Sam and Snowden gave a category theoretic approach to representation stability. They defined the stable representation category, $\text{Rep}^{st} (G_{*})$, and a specialization functor,
$$
\Gamma_{d}: \text{Rep}^{st} (G_{*}) \rightarrow \text{Rep} (G_{d}),
$$ 
mapping the stable representation category to the representation category for some group, $\text{Rep} (G_{d})$.

Their approach first examined the structure of the stable representation category and then examined the structure of the specialization functor. They did this for five families of groups: (1) the general linear groups, $GL_{n}$, (2) orthogonal groups, $O_{n}$, (3) symplectic groups, $Sp_{n}$, (4) symmetric group $S_{n}$, and (5) the general affine group, $GA_{n}$. They noted that although the results for each group are very similar, they have no unified theory.

This paper also laid the foundation for future combinatorial analysis of the categories Sam and Snowden established. In particular, the authors used relationships established via Schur-Weyl Duality to study $GL_{n}$, $O_{n}$, and $Sp_{n}$ via Brauer algebras and symmetric groups vis partition algebras.This approach echoed that of Bowman, De Visscher, and Orellana in \cite{BDO}, who used the stability of the partition algebra to study the Kronecker coefficients. A major result of \cite{SSc} was the use of the duality between the symmetric group and partition algebra to describe the infinite symmetric group, and also establishing the upwards and downwards partition categories, denoted (up) and (dp) respectively, to aid in their analysis. This led to one of the paper's most important results. Sam and Snowden define (up) and (dp) as the following.

\begin{defn}\cite[Section 6.]{SSc} (6.3.8) The \emph{downwards partition category} has objects
that are finite sets. A morphism $L \rightarrow L'$ is a partition of $L \amalg L'$ in which each part meets $L$. Given a
morphism $L \rightarrow L'$ represented by $\cU$ and a morphism $L' \rightarrow L''$  represented by $\cU ''$, the composition $L \rightarrow L''$  is represented by the partition obtained by gluing $\cU$ and $\cU ''$ along $L'$ and merging parts which meet. The \emph{upwards partition algebra} has everything reversed. 
\end{defn}

In \cite{SSc}[$\S$ 2.1.2], the authors let $\cA$ be an abelian category and the representation of $\Lambda$ is the functor $\Lambda \rightarrow \cA$. Let $\cA^{\Lambda}$ be the category of representations. The kernel, $\cK$ is an object of $\cA^{\Lambda}$. The contravariant functors $\Phi$ and $\Psi$ are given in \cite{SSc}[$\S$ 2.1.10] by
$$
\Phi : \text{Mod}^{f}_{\Lambda} \rightarrow \cA ,~~ \Phi(M) = \text{Hom}_{\Lambda}(M,\cK)
\quad \text{ and }\quad 
\Psi : \cA \rightarrow \text{Mod}_{\Lambda} ,~~ \Psi(N) = \text{Hom}_{\cA}(N,\cK).
$$
Define $\text{Mod}_{\text{(dp)}}^{\emph{f}}$ as the objects of  $\text{Mod}_{\text{(dp)}}$ that are of finite length. Let Rep$(S_{n})$ be the algebraic representations of the symmetric group and let
$$
\Phi : \text{Mod}_{\text{(dp)}}^{\emph{f}} \rightarrow \text{Mod}_{S_{n}}, \qquad \text{and} \qquad \Psi : \text{Mod}_{S_{n}} \rightarrow \text{Mod}_{\text{(dp)}}.
$$

\begin{thm} {\cite [Thm 6.3.30] {SSc}}
The functors $\Phi$ and $\Psi$ induce mutually quasi-inverse contravariant equivalences of tensor categories between $\mathrm{Rep(S_{n})}$ and $\mathrm{Mod^f_{(dp)}}$.
\end{thm} 

This theorem shows that the categories $\mathrm{Rep(S_{n})}$ and $\text{Mod}^f_{(dp)}$ are essentially equivalent, thereby establishing the relation between the category of representations of $S_n$ and the finitely generated modules of downward-partition algebra. 

In \cite{SSd}, the authors employed a combinatorial approach to studying certain properties of representations, such as Noetherianity, via combinatorial data. This work discussed a variety of interesting phenomena. In particular, we highlight results concerning TCAs in positive characteristic and analogues of the category FI.

\begin{defn} {\cite{SSd}}
Let $FI_{d}$ be the category whose objects are finite sets. Given two finite sets $S$ and $T$, a morphism $S \rightarrow T$ is a pair ($f$, $g$) where $f : S \rightarrow T$ is an injection and $g : T \ f(S) \rightarrow \{1, . . . , d\}$ is a
function.
\end{defn}

The authors were able to apply the results to several groups which served to illuminate the structure of those groups. Two such theorems are listed below but the authors are able to prove several more.

\begin{thm} {\cite [Cor 7.1.3] {SSd}}
If $R$ is left-Noetherian then $\text{Rep}_{R}(FI_{d})$ is Noetherian.
\end{thm}

\begin{thm} {\cite [Cor 7.3.6] {SSd}}
A TCA over a noetherian ring $R$ finitely generated in degree 1 is Noetherian.
\end{thm}

\pagebreak
\appendix

\section{List of finitely generated $\FI$-modules} \label{app:A.1}
The following list is a continuation of the list provided in \cite{CEFa}, extended to include examples developed in the growing body of literature.

\begin{center}
\begin{tabular}{ c c }  
 \textbf{$\FI$-module}: & \textbf{Description} \\ 
 \\
 $H^i(\mathrm{Conf}_n(M);\mathbb{Q})$ &  $\mathrm{Conf}_n(M)$ = configuration space of $n$ distinct ordered \\ & points on a connected, oriented manifold $M$ (\cite{CEFa}).\\
 \\
$R^{(r)}_J(n)$ & $J = (j_1,...j_r), R^{(r)}(n) = \bigoplus_J R_J^{(r)}(n)=r$-diagonal coinvariant algebra\\ & on $r$ sets of $n$ variables (\cite{CEFa}).\\
\\
$H^i(\mathcal{M}_{g,n}; \QQ)$ & $\mathcal{M}_{g,n}$ = moduli spaces of $n$-pointed genus with $g \geq 2$ curves (\cite{CEFa}). \\
\\
$R^i(\mathcal{M}_{g,n})$ & $i^{th}$ graded piece of tautological ring of $\mathcal{M}_{g,n}$ (\cite{CEFa}).\\
\\
$\mathcal{O}(X_{P,r}(n))_i$ & space of degree $i$ polynomials on $X_{P,r}(n)$, the rank variety\\ & of $n \times n$ matrices of $P$-rank $\leq r$ (\cite{CEFa}).\\
\\
$G(A_n/ \QQ)_i$ & degree $i$ part of the Bhargava-Satriano Galois closure of \\ & $A_n = \QQ[x_1,...,x_n]/(x_1,...,x_n)^2$ (\cite{CEFa}).\\
\\
$H^i(\mathcal{I}_n; \QQ)_{alb}$ & degree $i$ part of the subalgebra of $H^*(\mathcal{I}_n; \QQ)$ generated by $H^1(\mathcal{I}_n; \QQ)$,\\ & where $\mathcal{I}_n$ = genus $n$ Torelli group (\cite{CEFa}).
\\
\\
$H^i(IA_n; \QQ)_{alb}$ &  degree $i$ part of the subalgebra of $H^*(IA_n; \QQ)$ generated by $H^1(IA_n; \QQ)$,\\ & where $IA_n$ = Torelli subgroup of $\mathrm{Aut}(F_n)$ (\cite{CEFa}).\\
\\
$\mathrm{gr}(\Gamma_n)_i$ & $i^{th}$ graded piece of associated graded Lie algebra of groups $\Gamma_n$,\\ & including $\mathcal{I}_n, IA_n$, and pure braid group $P_n$ (\cite{CEFa}). \\
\\
$H^i(P\Sigma_n; \QQ)$ & $P\Sigma_n$ = pure string motion groups (\cite{CEFa}). \\
\\
$H^i(\mathcal{M}_W(n); \QQ)$ & $\mathcal{M}_W(n)$ = complexified complements of the reflecting hyperplanes $W_n$,\\ &
 where $W$ is either $S_n, BC_n$, or $D_n$.(\cite{W}) \\
\\
$H^i(\mathrm{Conf}_n(M);\mathbb{K})$ & where $\mathbb{K}$ is a Noetherian ring, and $M$ is a connected orientable manifold\\ & of dimension $\geq 2$ with the homotopy type of a finite $CW$ complex.(\cite{W})\\
\\
$E^{p,q}_*$ & Leray, Leray-Serre or Hochschild-Serre spectral sequence (\cite{JR}).\\
\\
$H^i(E; \mathbb{K})$ & $E$ is one of the above spectral sequences (\cite{JR}).\\
\\
$H^i(P\mathrm{Mod}_{g,r}^n; \QQ)$ & $P\mathrm{Mod}_{g,r}^n$ = pure mapping class group (\cite{JR}).\\
\\
$H^i(B \mathrm{ PDiff}^n(M); \QQ)$ & $B \mathrm{ PDiff}^n(M)$ = the classifying space of
the pure mapping class group (\cite{JR}).
\end{tabular}
\end{center}

\pagebreak


\begin{thebibliography}{1}

  \bibitem[BoDeOr]{BDO} C.\ Bowman, M.\ De Visscher, R.\ Orellana, \emph{The partition algebra and the Kronecker coeffiecients}, Trans. Amer. Math. Soc. 367 (2015), 3647-3667. arXiv:1210.5579 
  
  \bibitem[BR] {BR} H. \ Barcelo, A. \ Ram, \emph{Combinatorial representation theory}, (1997), 	arXiv: math/9707221

  \bibitem[Chu]{C} T.\ Church, \emph{Homological stability for configuration spaces of manifolds}, Inventiones mathematicae {\bf 188 No. 2} (2012), 465-504. arXiv:1103.2441

  \bibitem[ChEl]{CE} T.\ Church, J.\ Ellenberg, \emph{Homology of FI-modules}, (2014), arXiv: 1506.01022

  \bibitem[ChElFa1]{CEFa} T.\ Church, J.\ S.\ Ellenberg, B.\ Farb, \emph{FI-modules:  a new approach to stability for $S_n$-representations}, (2012). arXiv:1204.4533

  \bibitem[ChElFa2]{CEFb} T.\ Church, J.\ S.\ Ellenberg, B.\ Farb, \emph{FI-modules and stability for representations of symmetric groups}, Duke Math.\ J.\ {\bf 164 No. 9 } (2015), 1833-1910. arXiv:1204.4533v4

  \bibitem[ChElFaNa]{CEFN} T.\ Church, J.\ S.\ Ellenberg, B.\ Farb, R.\ Nagpal, \emph{FI-modules over Noetherian rings}, Geometry and Topology {\bf 18} (2014), 2951–2984. arXiv:1210.1854v2

  \bibitem[ChFa]{CF} T.\ Church, B.\ Farb, \emph{Representation theory and homological stability,}, Adv.\ Math.\ (2013), 250-314. arXiv:1008.1368
  
  \bibitem[GoWa]{GW} R. \ Goodman, N.R.  \ Wallach \emph{Symmetry, Representations, and Invariants}, (2009), Print.

  \bibitem[Far]{F} B.\ Farb, \emph{Representation stability}, to appear in the Proceedings of the 2014 Seoul ICM, (2014). arXiv: 1404.4065

 \bibitem[HaRa]{HaRa} T.\ Halverson, A.\ Ram, \emph{Partition Algebras}, (2004) arXiv:math/0401314 [math.RT]
  
  \bibitem[HeRe]{HR} P.\ Hersh, Vi.\ Reiner, \emph{Representation stability for cohomology of configuration spaces $\mathbb{R^d}$}, 	arXiv:1505.04196 [math.CO], 1 July 2015
  
  \bibitem[Rol]{JR} R.\ Jim\'{e}nez Rolland, \emph{On the cohomology of pure mapping class groups as FI-modules}, arXiv:1207.6828v2 [math.GT], 30 Sep 2013
  
  \bibitem[L]{LL} L.\ Li, \emph{Homological degrees of representations of categories with shift functors}, (2015). arXiv:1507.08023
  
  \bibitem[Mur1]{Ma} F.\ Murnaghan, \emph{The analysis of the Kronecker product of irreducible representations
of the symmetric group}, Amer.\ J.\ Math.\ {\bf 188 No. 60} (1938), no.3 761-784.

  \bibitem[Mur2]{Mb} F.\ Murnaghan, \emph{On the analysis of the Kronecker product of irreducible representations of
$S_{n}$}, Proc.\ Nat.\ Acad.\ Sci.\ U.S.A.\ {\bf 188 No. 41} (1955), no.3 515-518.

\bibitem[Ra] {R} E. \ Ramos, \emph{Homological invariants of $\FI$-modules and $\FI_G$-modules}, (2015). arXiv:1511.03964

\bibitem[Nag]{N} R.\ Nagpal, \emph{FI-modules and the cohomology of modular representations of symmetric groups}, arXiv:1505.04294
  
  \bibitem[Put]{P} A.\ Putnam \emph{Stability in the homology of congruence subgroups}, (2015). arXiv:1201.4876

  \bibitem[PuSa]{PS} A.\ Putnam, S.\ V.\ Sam, \emph{Representation stability and finite linear groups}, (2014). arXiv:1408.3694
 

  \bibitem[SaSn1]{SSa} S.\ Sam, A.\ Snowden, \emph{GL-equivariant modules over polynomial rings in infinitely many variables}, (2012), arXiv: 1206.2233

  \bibitem[SaSn2]{SSb} S.\ Sam, A.\ Snowden, \emph{Introduction to Twisted Commutative Algebras}, (2012), arXiv:1209.5122 

  \bibitem[SaSn5]{SSf} S.\ Sam, A.\ Snowden, \emph{Proof of Stembridge's conjecture on stability of Kronecker coefficients}, Journal of Algebraic Combinatorics, to appear, arXiv:1501.00333v2

  \bibitem[SaSn3]{SSc} S.\ Sam, A.\ Snowden, \emph{Stability patterns in representation theory}, (2013), arXiv:1302.5859

  \bibitem[SaSn4]{SSd} S.\ Sam, A.\ Snowden, \emph{Gr$\ddot{o}$bner methods for representations of combinatorial categories}, (2014), arXiv:1302.5859 

  \bibitem[NaSaSn]{SSe} R. Nagpal, S.\ Sam, A.\ Snowden, \emph{Noetherianity of some degree two twisted commutative algebras}, (2015), arXiv:1302.5859 

  \bibitem[Wil1]{Wil} J.\ C.\ H.\ Wilson , {\em FI$_{\emph{W}}$ modules and stability criteria for representations of Weyl groups}, J.\ Algebra {\bf 420} (2014), 269-332. arXiv: 1309.3817

  \bibitem[Wil2]{W} J.\ C.\ H.\ Wilson, \emph{$FI_W$-modules and constraints on classical Weyl group characters}, arXiv:1503.08510 [math.RT]

  \end{thebibliography}
\end{document}